\documentclass[runningheads]{llncs}
\usepackage{graphicx}
\usepackage{amsfonts,amsmath}
\newcommand{\RR}{{\mathbb R}}

\begin{document}

\title{Crossings Between Non-homotopic Edges
\thanks{Supported by the National Research, Development and Innovation Office, 
NKFIH, KKP-133864, K-131529, K-116769, K-132696, 
by the Higher Educational Institutional Excellence Program 2019 
NKFIH-1158-6/2019, 
the Austrian Science Fund (FWF), grant Z 342-N31, 
by the Ministry of Education and Science of the Russian Federation 
MegaGrant No.\ 075-15-2019-1926, 
and by the ERC Synergy Grant ``Dynasnet'' No.\ 810115.}}
%
%
\author{J\'anos Pach\inst{1,2,4}
\and
G\'abor Tardos\inst{1,3,4}
\and
G\'eza T\'oth\inst{1,5}
}
\authorrunning{J. Pach et al.}
%
\institute{
R\'enyi Institute, Budapest\email{\{pach, tardos, geza\}@renyi.hu} \and
IST Austria, Vienna \and
Department of Mathematics, Central European University, Budapest \and
Moscow Institute of Physics and Technology, Moscow \and
Budapest University of Technology and Economics, SZIT, Budapest}
\maketitle              
\begin{abstract}
We call a multigraph {\em non-homotopic} if it can be drawn in the plane
  in such a way that no two edges connecting the same pair of vertices can
  be continuously transformed into each other without passing through a
  vertex, and no loop can be shrunk to its end-vertex in the same way.
  It is easy to see that a non-homotopic multigraph on $n>1$ vertices
  can have arbitrarily many edges. We prove that the number of crossings between
  the edges of a non-homotopic multigraph with $n$ vertices and $m>4n$ edges
  is larger than $c\frac{m^2}{n}$ for some constant $c>0$, and that this bound
  is tight up to a polylogarithmic factor. We also show that the lower bound is
  not asymptotically sharp as $n$ is fixed and $m\longrightarrow\infty$
\end{abstract}

\section{Introduction}\label{intro}

A standard parameter for measuring the non-planarity of a graph $G$ is its {\em crossing number},
which is defined as the smallest number ${\rm cr}(G)$ of crossing points in any drawing of $G$
in the plane.
For many interesting variants of the crossing number, see  ~\cite{Sch13,Sch18,Sz04,PT00b}.
Computing ${\rm cr}(G)$ is an NP-complete problem~\cite{GJ83}.

Perhaps the most useful result on crossing numbers, is the
so-called {\em crossing lemma},
proved independently by Ajtai, Chv\'atal, Newborn,
Szemer\'edi~\cite{ACNS82} and Leighton~\cite{L83},
according to which the crossing number of any graph
with $n$ vertices and $m>4n$ edges is at least
$c\frac{m^3}{n^2}$, for a suitable constant $c>0$.
For the best known value of the constant $c$, see \cite{PRTT06,A19}.
This result, which is tight up to the constant factor, has been
successfully applied to a variety of problems in discrete
and computational geometry, additive
number theory, algebra, and elsewhere~\cite{D98,Sz97}.
In some applications, it was the bottleneck that one needed a lower bound
on the crossing number of a {\em multigraph} rather than a graph. Obviously,
the crossing lemma does not hold in this case, as stated. Indeed, one can
connect a pair of vertices ($n=2$) with $m$ parallel edges without creating
any crossing. However, for multigraphs $G$ with maximum edge {\em multiplicity}
$k$ and $m>4kn$ edges, Sz\'ekely~\cite{Sz97} established the lower bound
${\rm cr}(G)> c'\frac{m^3}{kn^2}$, where $c'>0$ is another constant.
This bound is also tight, up to the constant factor.
\'Agoston and P\'alv\"olgyi \cite{AP20} observed that $c'$ can
be chosen to be the same as the best known constant $c$ in the
crossing lemma (presently, $\frac{1}{29}$).

As the multiplicity $k$ increases, Sz\'ekely's bound gets weaker and weaker.
Luckily, the term $k$ in the denominator can be eliminated in several
special cases; see \cite{PT18,KPTU18}. That is, the result holds without
putting any upper bound on the edge multiplicity. However, in all of these cases, we have
to assume (among other things) that no two adjacent edges cross.
\smallskip

In this paper, we study the analogous question under the
weakest possible assumption.
Obviously, we need to assume that no pair of parallel
edges or loops are
{\em homotopic}, i.e., they cannot be continuously
deformed into each
other so that their interiors do not pass through
any vertex.
As we have noted above above, without this assumption, a
multigraph can have arbitrarily many non-crossing edges.
For simplicity, we will also assume that there are no
{\em trivial} loops, that is, no loop can be transformed
into a point. Clearly, this latter assumption can be eliminated
as the first condition already implies that there is at
most a single trivial loop at any vertex.

\smallskip

To state our results, we need to agree about the definitions.

A {\em multigraph} is a graph in which parallel edges and loops are permitted.
A {\em topological graph (or multigraph)} is a graph (multigraph)
$G=(V,E)$ drawn in the plane with the property that every vertex is
represented by a distinct point and every edge $e\in E$ is represented by a
continuous curve, i.e., a continuous function $f_e\colon[0,1]\to \RR^2$ with
$f_e(0)$ and $f_e(1)$ being the endpoints of $e$. In terminology,
we do not distinguish between the vertices and the points representing them.
In the same spirit, if there is no danger of confusion, we often use the term edge
instead of the curve $f_e$ representing it or the image of $f_e$.
As we deal with non-oriented multigraphs, we treat the functions $f_e(t)$ and $f_e(1-t)$ as being the same.
We assume that no edge passes through any vertex (i.e., $f_e(t)\not\in V$ for $0<t<1$).

\smallskip

The {\em crossing number} of a {\em topological multigraph} $G$ is the number of crossings between
its edges, i.e, the number of unordered pairs of distinct pairs $(e,t),(e',t')\in E\times(0,1)$
with $f_{e}(t)=f_{e'}(t')$. With a slight abuse of notation,
this number will be denoted also by ${\rm cr}(G)$.

Two parallel edges, $e, e',$ connecting the same pair of vertices, $u,v\in V$ are
{\em homotopic}, if there exists a continuous function ({\em homotopy}) $g\colon[0,1]^2\to\RR^2$
satisfying the following three conditions.
     $$g(0,t)=f_{e}(t) \mbox{ and } g(1,t)=f_{e'}(t) \mbox{ for all } t\in[0,1],$$
     $$g(s,0)=u \mbox{ and } g(s,1)=v \mbox{ for all } s\in[0,1],$$
     $$g(s,t)\not\in V \mbox{ for all } s,t\in(0,1).$$
Recall that we do not distinguish $f_{e}$ from $f_{e}(1-t)$, so we call $e$ and $e'$ homotopic also if $f_{e}(1-t)$ and $f_{e'}(t)$ are homotopic in the above sense. A loop at vertex $u$ is said to be \emph{trivial} if it is homotopic to the constant function $f(t)=u$.

\smallskip

A topological multigraph  $G=(V,E)$ is called {\em non-homotopic} if
it does not contain two homotopic edges, and does not contain a trivial loop.

Obviously, if $G$ is a simple topological graph (no parallel edges or loops), then it is non-homotopic.
A non-homotopic multigraph with zero or one vertex has no edge. However, if the number of vertices $n$ is at least $2$, the number of edges can be arbitrarily large, even infinite. Our first result provides a lower bound on the crossing number
of non-homotopic topological multigraphs in terms of the number of their vertices and edges.

\begin{theorem}\label{lowerbound1}
The crossing number of a non-homotopic topological multigraph $G$ with $n>1$ vertices and $m>4n$ edges satisfies ${\rm cr}(G)\ge\frac1{24}\frac{m^2}{n}$.
\end{theorem}

This bound is tight up to a polylogarithmic factor.

\begin{theorem}\label{upperbound}
  For any
$n\ge 2$, $m>4n$,
there exists a non-homotopic multigraph
  $G$ with $n$ vertices and $m$ edges
  such that its crossing number satisfies
${\rm cr}(G)\le30\frac{m^2}{n}\log_2^2\frac mn$.
\end{theorem}

\noindent The constant $30$ in the theorem was chosen for the proof to work for all $n$ and $m$,
and we made no attempt to optimize it. However, it can be replaced by $1+o(1)$ if both $n$ and $m/n$ go to infinity.
\smallskip

Define the function ${\rm cr}(n,m)$ as the minimum crossing number
of a non-homotopic multigraph with $n$ vertices and $m$ edges.
Theorems~\ref{lowerbound1} and \ref{upperbound} can be stated as
$$\frac1{24}\frac{m^2}{n}\le{\rm cr}(n,m)\le 30\frac{m^2}n\log_2^2\frac mn,$$
for any $n\ge2$ and $m>4n$.
We have been unable to close the gap between these bounds.
Our next theorem shows that the lower bound is not tight.

\begin{theorem}\label{lowerbound2}
The minimum crossing number of a non-homotopic multigraph with
$n\ge 2$ vertices and $m$ edges is super-quadratic in $m$.
That is, for any fixed $n\ge 2$, we have
 $$\lim_{m\rightarrow\infty}\frac{{\rm cr}(n,m)}{m^2}=\infty.$$
More precisely, we obtain
$\frac{{\rm cr}(n,m)}{m^2}=\Omega({\log m}^{1/(6n)}/{n^7})$.
\end{theorem}


Let $n, k$ be positive integers, and consider a set $S$ obtained
from the Euclidean plane by removing $n$ distinct points.
Fix a point $x\in S$. An oriented loop in $S$ that starts and
ends at $x$ is called an
{\em $x$-loop}. An $x$-loop may have self-intersections. Contrary to our convention for edges of a topological multigraph, we do distinguish between an $x$-loop and its reverse. We consider the {homotopy type of $x$-loops} in $S$, that is, we consider two loops {\em homotopic} if one can be continuously transformed to the other within $S$. When counting self-intersections of $x$-loops or intersections between two $x$-loops, we count points of multiple intersections with the appropriate multiplicity.
\smallskip

To establish Theorems \ref{upperbound} and \ref{lowerbound2},
we study the following topological problem of independent
interest \cite{JMM96}.

\begin{problem}
Let $n,k\ge 1$ be integers, let $S$ denote the set obtained from $\RR^2$ by removing $n$ distinct points, and let us fix $x\in S$. Determine or estimate the maximum number $f(n,k)$ of pairwise non-homotopic $x$-loops in $S$ such that none of them passes through $x$, each of them has fewer than $k$
self-intersections and every pair of them cross fewer than $k$ times.
\end{problem}

It is not at all obvious that $f(n,k)$ is finite. However,
in the sequel we show that this is the case.
This fact is crucially important for the proof of Theorem~\ref{lowerbound2}.

\begin{theorem}\label{diffloop}
For any integers $n\ge2$ and $k\ge1$, we have
$$f(n,k)<2^{(2k)^{2n}}.$$
\end{theorem}

Our proof of Theorem~\ref{upperbound} is based on a lower bound on $f(n,k)$. For this application, all we need is the $n=2$ special case. Next we state a lower bound valid for all $n$. 

\begin{theorem}\label{diffloop2}
Let $n\ge 2$ and $k\ge 1$ be integers. If $2\le n\le2k$, then $$f(n,k)\ge2^{\sqrt{nk}/3}$$ holds. For $n\ge2k$, we have
$$f(n,k)\ge(n/k)^{k-1}.$$
\end{theorem}

There is a huge gap between this bound and the upper bound in Theorem~\ref{diffloop}.
We suspect that the truth is
to the lower bound. More precisely, we conjecture that $\log f(n,k)$ can be
bounded from above by a polynomial of $k$ whose degree does not depend on $n$.
For $n=2$ we have $2^{\sqrt{k}/3}\le f(2, k)\le 2^{16k^4}$.

\smallskip


Our paper is organized as follows. In Section~2,
we establish Theorem~\ref{lowerbound1}. In Section~3, we present some constructions proving Theorem~\ref{diffloop2},
and apply them to deduce Theorem~\ref{upperbound}.
In Section 4, we prove Theorem \ref{diffloop}.
The proof of Theorem \ref{lowerbound2} is moved to the
Appendix. 


\section{Loose Multigraphs---Proof of Theorem~\ref{lowerbound1} }\label{mainbounds}


One can also define topological multigraphs and non-homotopic
multigraphs on the sphere $S^2$. If we consider $S^2$ as the
single point compactification of the plane with the \emph{ideal point}
$p^*$, then any topological multigraph $H$ drawn in the plane remains a
topological multigraph on the sphere. However, it may lose the
non-homotopic property, as the addition of the ideal point
$p^*$ may turn a loop trivial or two parallel edges homotopic.
This can be avoided by adding $p^*$ as an isolated vertex to $H$:
in this case, the resulting multigraph $H^*$ is non-homotopic even on the sphere.

\vskip0.5cm

\begin{figure}[!ht]
\begin{center}
\scalebox{0.5}{\includegraphics{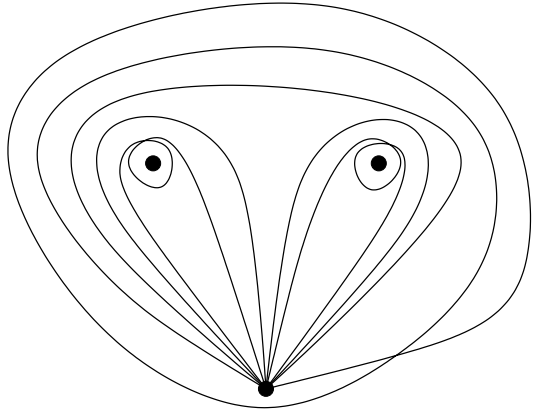}}
\caption{A non-homotopic loose multigraph with $3$ vertices and $6$
  edges (loops).
}\label{6loop}
\end{center}
\end{figure}

\smallskip

We say that a topological multigraph is {\em loose} if no pair of distinct
edges cross each other. An edge (in particular, a loop) is allowed
to cross itself. We start by finding the maximum number of edges in
a loose non-homotopic multigraph on the sphere or in the plane, for a given number of vertices.
We will see that despite allowing parallel edges, loops,
and self-intersections, loose non-homotopic multigraphs
with $n>2$ vertices on the sphere cannot have more than $3n-6$ edges,
the maximum number of edges of a simple planar graph.
However, there are many other nontrivial examples,
for which this bound is tight.
The interested reader can verify that, for all
$n>2$, there are extremal examples, all of whose edges are loops. See Fig.~1 for the case of three vertices in the plane.

\begin{lemma}\label{claim0}
On the sphere, any loose non-homotopic multigraph with $n>2$ vertices has $m\le3n-6$ edges.
For $n=2$, the maximum number of edges is 1.
\end{lemma}


The proof of Lemma \ref{claim0} can be found in the Appendix.

\begin{lemma}\label{claim1}
In the plane, any loose non-homotopic multigraph with $n\ge1$ vertices has at most $3n-3$ edges. This bound can be achieved for every $n$.
\end{lemma}

\noindent{\bf Proof.} Let $H$ be a loose non-homotopic multigraph in the plane with $n\ge1$ vertices and $m$ edges. Consider the plane as the sphere $S^2$ with a point $p^*$ removed. Add $p^*$ to $H$ as an isolated vertex, to obtain a topological multigraph $H'$ on the sphere. Then $H'$ is a loose non-homotopic multigraph with $n+1$ vertices and $m$ edges. If $n>1$, applying Lemma~\ref{claim0} to $H'$, we obtain that $m\le3n-3$, as required. If $n=1$, then $H$ is a single-vertex topological multigraph in the plane, so all of its edges must be trivial loops. However, by definition, a non-homotopic multigraph cannot have any trivial loop. This completes the proof of the upper bound.
\smallskip

There are many different constructions for loose non-homotopic multigraphs for which the bound in the lemma is achieved. Such a topological multigraph may have several components and several self-intersecting loops. (However, all self-crossings of non-loop edges must be ``homotopically trivial'':  the removal of the closed curve produced by such a self-crossing does not change the homotopy type of the edge.)
\smallskip

Here, we give a very simple construction. If $n>2$, we start with a triangulation with $n$ vertices and $3n-6$ edges. Let $uvw$ be the boundary of the unbounded face. Add another non-self-intersecting edge connecting $u$ and $v$ in the unbounded face, which is not homotopic with the arc $uv$ of $uvw$. Finally, we add two further loops at $u$. First, a simple loop $l$ that has all other edges and vertices (except $u$) in its interior, and then another loop $l'$  outside of $l$, which goes twice around $l$. (Of course, $l'$ must be self-intersecting.)

If $n=1$, the graph with no edge achieves the bound of the lemma. For $n=2$, draw an edge $e$ connecting the two vertices, $u$ and $v$. Then add two loops at $u$, as above: a simple loop $l$  around $e$ and another loop $l'$ that winds around $l$ twice. \hfill  $\Box$
\medskip

\noindent{\bf Proof of Theorem~\ref{lowerbound1}.} Let $G$ be a non-homotopic topological
multigraph in the plane with $n>1$ vertices and $m>4n$ edges.

Let $D$ denote the {\em non-crossing graph} of the edges of $G$, that is, let $V(D)=E(G)$ and
connect two vertices of $D$ by an edge if and only if the corresponding edges of $G$
do not share an interior point. Any clique in $D$ corresponds to a loose non-homotopic sub-multigraph of $G$. Therefore, by Lemma~\ref{claim1}, $D$ has no clique of size $3n-2$.
Thus, by Tur\'an's theorem~\cite{Tu41},
$$|E(D)|\le\frac{|V(D)|^2}{2}\left(1-\frac1{3n-3}\right)=\frac{m^2}{2}\left(1-\frac1{3n-3}\right).$$
The crossing number ${\rm cr}(G)$ is at least the number of crossing pairs of edges in $G$, which is equal to the number of non-edges of $D$. Since $m>4n$, we have
$${\rm cr G}\ge\binom m2-\frac{m^2}{2}\left(1-\frac1{3n-3}\right)\ge\frac1{24}\frac{m^2}{n},$$
as claimed. \hfill   $\Box$
\medskip

The proof above gives a lower bound on the number
of crossing pairs of edges in $G$,
and in this respect it is tight up to a constant factor.
To see this, suppose for simplicity that $n$ is even and $m$ is divisible by $n$.
Let $G_0$ be a non-homotopic topological multigraph with
two vertices and $\frac{2m}n$ non-homotopic loops on one of its vertices.
Taking $\frac{n}2$ disjoint copies of
$G_0$, we obtain a non-homotopic topological
multigraph with $n$ vertices, $m$ edges, and $<\frac{m^2}{n}$
crossing pairs.


\section{Two Constructions---Proofs of Theorems~\ref{diffloop2} and \ref{upperbound}}

The aim of this section is to demonstrate how to construct topological graphs with many edges and families consisting  of many loops, without creating many crossings. The constructions are based on the description of the fundamental group of the plane from which a certain number of points have been removed.
\medskip

\noindent{\bf Proof of Theorem~\ref{diffloop2}.}
Let $S=\RR^2\setminus\{a_1,\dots,a_n\}$, where $a_1,\ldots, a_n$ are distinct points in the plane, and let $x\in S$ be also fixed. Assume without loss of generality that $a_i=(i, 0),\; 1\le i\le n,$ and $x=(0, -1)$. Recall that an {\em $x$-loop} is a (possibly self-crossing)
oriented path in $S$ from $x$ to $x$, i.e., a continuous function
$f\colon[0,1]\to S$ with $f(0)=f(1)=x$.

Note that the homotopy group of $S$ is the free group $F_n$ generated by $g_1,\dots,g_n$, where $g_i$ can be represented by a triangular $x$-loop around $a_i$, for example the one going from $x$ to $(2i-1, 1)$, from here to $(2i+1, 1)$, and then back to $x$ along three straight-line segments; see~\cite{Ly77}.


\smallskip

We define an {\em elementary loop} to be a polygonal $x$-loop with intermediate vertices $$(1,\pm1/2),\,(2, \pm1/2),\,\dots,\,(n,\pm1/2),\,(n+1,-1),$$ in this order. There are $2^n$ distinct elementary loops, depending on the choice of the signs. Each of them represents a distinct homotopy class of the form $g_{i_1}\cdots g_{i_t}$, where the indices form a strictly increasing sequence. By making infinitesimal perturbations on the interior vertices of the elementary loops, we can make sure that every pair of them intersect in at most $n-1$ points. Thus, we have $f(n,n)\ge2^n$.
\smallskip

We call $1\le i<n$ a \emph{sign change} in the elementary loop $l$ if $l$ passes through both $(i,1/2)$ and $(i+1,-1/2)$, or both $(i,-1/2)$ and $(i+1,1/2)$. There are precisely $2\binom{n-1}j$ elementary loops with exactly $j$ sign changes. The reader can easily verify that crossings between perturbed elementary loops are unavoidable only if a sign change occurs. More precisely, for $k\le n$, one can perturb all elementary loops with at most $k-1$ sign changes such that every pair cross at most $k-1$ times. Hence, we have $f(n,k)\ge2\sum_{j=0}^{k-1}\binom{n-1}j\ge2\binom n{k-1}>(n/k)^{k-1}$, completing the proof of the theorem, whenever $n\ge2k$.
\smallskip

If $k\le n\le2k$, we have $f(n,k)\ge f(k,k)\ge2^k<2^{\sqrt{nk}/3}$. Similarly, if $n\le k\le9n$, we have $f(n,k)\ge f(n,n)\ge 2^n\ge2^{\sqrt{nk}/3}$, and we are done.
\smallskip

Finally, in the case $k>9n$, we consider all $x$-loops which can be obtained as the product (concatenation) of $j=\lfloor\sqrt{\frac{k-1}n}\rfloor\ge3$ elementary loops. Unfortunately, some of these concatenated $x$-loops will be homotopic. For example, if the elementary loops $l_1,l_2,l_3,$ and $l_4$ represent the homotopy classes $g_1,g_2g_3,g_1g_2,$ and $g_3$, respectively, then $l_1l_2$ and $l_3l_4$ are homotopic. To avoid this complication, we only use the $2^{n-1}$ elementary loops that represent homotopy classes involving $g_1$ (that is, the ones with $(1,+1/2)$ as their first intermediate vertex). Then no two of the resulting $2^{j(n-1)}$ $x$-loops will be homotopic. By infinitesimal perturbation of the interior vertices of these $x$-loops (including the $j-1$ interior vertices at $x$), we can attain that they do not pass through $x$, and no two polygonal paths corresponding to a single elementary loop intersect more than $n$ times. Therefore, any pair of perturbed concatenated loops cross at most $j^2n<k$ times, and the same bound holds for the number of self-intersections of any concatenated loop.
This yields that $f(n,k)\ge2^{j(n-1)}\ge 2^{\sqrt{nk}/3}$.  \hfill$\Box$
\medskip

\noindent{\bf Proof of Theorem~\ref{upperbound}.}
We want to construct a non-homotopic topological multigraph $G$ with $n$ vertices, $m$ edges, and few crossings. We distinguish 3 cases.
\smallskip

{\em Case A:} If $n=3$, we set $k=\lceil2\log_2^2(2m)\rceil$. Theorem~\ref{diffloop2}  guarantees that $f(2,k)\ge2m$. Thus, there are $2m$ pairwise non-homotopic $x$-loops in $S=\RR^2\setminus\{a_1,a_2\}$ such that each of them has fewer than $k$ self-intersections and any pair intersect fewer than $k$ times. Regard this arrangement as a topological multigraph $G$ with $2m$ edges on the vertex set $\{a_1,a_2,x\}.$ All edges are $x$-loops. At most one of them is trivial, and for each loop edge there is at most one other loop edge homotopic to it (which must come from an $x$-loop with inverse orientation). Therefore, we can always select $m$ edges that form a non-homotopic multigraph. Obviously, we have ${\rm cr}(G)<k(m+\binom m2)$.
\smallskip

{\em Case B:} If $n>3$, we set $n^*=\lfloor n/3\rfloor$, $m_0=\lceil m/n^*\rceil$. Take $n^*$ disjoint copies of the non-homotopic multigraph $G_0$ with $3$ vertices and $m_0$ edges constructed in Case~A. We add at most $2$ isolated vertices and remove a few edges if necessary to obtain a non-homotopic multigraph on $n$ vertices and $m$ edges. We clearly have ${\rm cr}(G)\le n^*{\rm cr}(G_0)$.

Clearly, the crossing numbers of the graphs constructed in Cases~A and B are within the bound stated in the theorem.
\smallskip

{\em Case C:} If $n=2$, we cannot use Theorem~\ref{diffloop2} directly. Note that all edges of the non-homotopic multigraphs $G$ constructed in Case~A were loops at a vertex $x$, and these $x$-loops were pairwise non-homotopic even in the set obtained from the plane by keeping $x$, but removing every other vertex. Now we cannot afford this luxury without creating $\Omega(m^3)$ crossings. However, even for $n=2$, we can construct a topological multigraph $G$ with many pairwise non-homotopic edges and relatively few crossings, as sketched below.
\smallskip

Let $V(G)=\{a_1,a_2\}$, where $a_1$ and $a_2$ are distinct points in the plane, and set $S=\RR\setminus V(G)$. Choose a base point $x\in S$ not on the line $a_1a_2$. Now the homotopy group of $S$ is the free group generated by two elements, $g_1$ and $g_2$, that can be represented by triangular $x$-loops around $a_1$ and $a_2$, respectively. By the proof of Theorem~\ref{diffloop}, with the notation used there, we can construct $2^j$ pairwise non-homotopic $x$-loops in $S$ with few crossings. Each of these $x$-loops, $l$, can be turned into either a loop edge at the vertex $a_1$ or into an $a_1a_2$ edge, as follows: we start with the straight-line segment $a_1x$, then follow $l$, finally add a straight-line segment from $x$ to either $a_1$ (for a loop edge) or to $a_2$ (to obtain a non-loop edge). After infinitesimally perturbing the resulting edges, one can easily bound the crossing number. However, now we face a new complication: there may be a large number of pairwise homotopic edges. In Case~A, when we regarded $x$-loops as loop edges in a topological multigraph having $x$ as a vertex, two loop edges could only be homotopic if the corresponding $x$-loops represented the same or inverse homotopy classes. Now the situation is more complicated: a loop edge constructed from an $x$-loop representing an element $g$ in the homotopy group is homotopic to an another edge constructed from another $x$-loop representing $g'$ if and only if we have $g'=g_1^sgg_1^t$ or $g'=g_1^sg^{-1}g_1^t$ for some integers $s$ and $t$. (For non-loop edges the corresponding condition is $g'=g_1^sgg_2^t$.)
We may have constructed more than two (even an unbounded number of) homotopic edges, but a closer look at
the $2^j$ $x$-loops constructed in the proof of the lower bound on $f(2,k)$ reveals that $2^{j-2}$
of them yield pairwise non-homotopic edges. \hfill$\Box$

\smallskip

\noindent{\bf Remark.} For $n\ge 3$, in our constructions all edges are loops.
By splitting the base points of the loops, we can
get constructions with no loops. 



\section{Loops with Bounded Number of Pairwise Intersections\\---Proof of Theorem~\ref{diffloop}}


Consider a loop (oriented closed curve) $l$ in the plane,
and a point $r$ not belonging to $l$.
The {\em winding number} of $l$ around $r$ is the number
of times the loop goes around $r$ in the counter-clockwise
direction. Going around $r$ in the clockwise direction counts negatively.

Let $S$ be obtained by removing a single point
$r$ from the plane. It is well known that two loops in
$S$ are homotopic if and only if their winding numbers around $r$ are the same.

\begin{lemma}\label{wind}
Let $l$ be any loop in the plane with fewer than $k$ self-intersections, and let $x$ be a point that does not belong to $l$.
Then the absolute value of the winding number of $l$ around $x$ is at most $k$.
\end{lemma}

\noindent{\bf Proof.}
Removing the image of $l$ from the plane, it falls into connected components,
called \emph{faces}. Obviously, the winding number of $l$ is the same around
any two points, $x$ and $y$, that belong to the same face.
Take a point in each face and connect two distinct points if the corresponding faces have a common boundary curve.
We get a connected graph.
If $x$ and $y$ are adjacent, then the winding number of $l$ around $x$ and $y$
differs by precisely $1$.
As $l$ has fewer than $k$ self-intersections, the number of faces is at most $k+1$.
The winding number of $l$ around any point of the unbounded face is zero.
Therefore, the winding
number of $l$ around any point not belonging to $l$ is between $-k$ and $+k$, as claimed. \hfill    $\Box$




\begin{corollary}\label{bazis}
For any integer $k>0$, we have $f(1,k)\le2k+1.$
\end{corollary}

\noindent{\bf Proof.} Let $x$ and $a$ be two distinct points in the plane $\RR^2$. Any two $x$-loops in $S=\RR^2\setminus\{a\}$ are homotopic if they have the same winding number around $a$. For an $x$-loop with fewer than $k$ self-crossing this is winding number is takes values between $-k$ and $k$. therefore, any collection of pairwise non-homotopic such loops has cardinality at most $2k+1$. \hfill $\Box$
\medskip

In the rest of this section, we estimate the function $f(n,k)$ for $n>1$. By the definition of $f(n,k)$, we have to consider a set $S$ that can be obtained from $\RR^2$ by removing $n$ distinct points.
As before, we consider the 2-sphere $S^2$ as the compactification of the plane with a single point $p^*$, the ``ideal point.''. To simplify the presentation, we view $S$ as a set obtained from $S^2$ by the removal of a set $T$ of $n+1$ points (including $p^*$). We also fix the common starting point $x\in S$ of all loops in $S$ that we consider.
\smallskip

Let $L$ be a collection of loops in $S$. The connected components of $S^2$ minus the set of all points of the elements of $L$ are called {\em $L$-faces}. Obviously, all $L$-faces are homeomorphic to the plane and the points of $T$ are scattered among them. We call $L$ {\em balanced} if no $L$-face contains $n$ or $n+1$ points of $T$.

\begin{lemma}\label{suit}
Let $k$ be a positive integer, let $x\in S$, and let $H$ be a collection of pairwise non-homotopic nontrivial $x$-loops in $S$, each of which has fewer than $k$ self-intersections.

If $|H|>2k+1$, then there is a balanced pair, $L$, of loops in $H$.
\end{lemma}

\noindent{\bf Proof.}
For $n=1$, there is no balanced family. Nevertheless, formally the statement holds even in this case, because Lemma~\ref{wind} implies that $|H|\le2k+1$. (In fact, now $|H|<2k+1$ because of the non-triviality condition.)
\smallskip

Suppose that $n>1$. Consider any loop $l\in H$. If $\{l\}$ is balanced, then any pair containing $l$ is also balanced and we are done. Otherwise, there is an $\{l\}$-face $F$ containing at least $n$ points of $T$. It cannot contain all points of $T$, because then $l$ would be contractible, that is, trivial.

Therefore, we can assume that there is a single point $t\in T$ outside $F$. We say that $l$ {\em separates} $t$. If two loops, $l_1, l_2\in H$, separate distinct points, $t_1, t_2\in T$, respectively, then $L=\{l_1,l_2\}$ is a balanced pair, because $t_1$ and $t_2$ must lie in separate $L$-faces,
distinct from all $L$-faces containing other points of $T$.
\smallskip

Hence, we may assume that all loops in $H$ separate the
same point $t\in T$. By symmetry, we may also assume $t\neq p^*$
($t$ is not the ideal point), so $t$ is in the plane. By Lemma~\ref{wind},
the winding number of any loop $l\in H$ around $t$ is between $-k$ and $+k$. If $|H|>2k+1$, by the pigeonhole principle, there are two distinct loops, $l_1,l_2 \in H$, with the same winding number around $t$. This implies that $L=\{l_1,l_2\}$ is a balanced pair. Indeed, otherwise all points in $T\setminus\{t\}$ would be in the same $L$-face $F$. In this case, all points of $T\setminus\{t,r\}$ would lie in the unbounded face of the arrangement of the loops $l_1$ and $l_2$ in the plane. Since $l_1$ and $l_2$ have the same winding number around $t$,
it would follow that they are homotopic, a contradiction.  \hfill $\Box$
\medskip

Now we are in a position to establish the following recurrence relation for $f(n,k)$,
which, together with Lemma~\ref{bazis}, implies the upper bound in Theorem~\ref{diffloop}.

\begin{lemma}\label{rekurzio}
For any $n>1$, $k>0$, we have
$f(n,k)\le(6kf(n-1,k))^{2k}$.
\end{lemma}

\noindent{\bf Proof.}
Consider a family $H$ of loops for which in the definition of $f(n,k)$ the maximum  is attained. That is, consider $S=S^2\setminus T$ with $|T|=n+1$, fix a point $x\in S$, and let $H$ consist of $f(n,k)$ pairwise non-homotopic $x$-loops in $S$ not passing through $x$, such that each loop has fewer than $k$ self-intersections and every pair of loops intersect in fewer than $k$ points. We may also assume, by infinitesimal perturbations, that there is no triple-intersection and that any intersection point of the loops is a {\em transversal} crossing, where one arc passes from one side of the other arc to the other side.
\smallskip

If $|H|\le2k+1$, we are done. Suppose that $|H|>2k+1$.
By Lemma~\ref{suit}, there exists a balanced two-element subset $L\subset H$.
Fix such a subset $L$,
and turn the arrangement of the two loops in $L$ to a
multigraph drawn on the sphere, as follows.
Regard $x$ and all intersection points as vertices,
so that we obtain a planar drawing of a 4-regular connected multigraph $G$ with at most $3k-2$ vertices.
Thus, the number of edges of $G$ satisfies $|E(G)|\le 6k-4$. For any edge $f$ of $G$, we designate an arbitrary curve in $S$ starting at $x$ and ending at an internal point of $f$, and we call it the {\em leash} of $f$. (For example, we may choose the leash to pass very close to the edges of $G$.) For any internal point $a$ of $f$, the {\em standard path from $x$ to $a$} is the leash of $f$ followed by the piece of $f$ between the endpoint of the leash and $a$. When referring to the standard path to $x$, we mean the single point curve.
\smallskip

Consider a loop $l\in H\setminus L$. Here $l$ starts at $x$ and has later some $j\le2k-2$ further intersections with the edges of $G$, each time properly crossing an edge from one face of $G$ to another. Define the {\em signature of $l$} as the sequence of these $j$ edges of $G$, together with the information where the loop starts and ends in a tiny neighborhood of $x$. For the latter, we just record the cyclic order of the initial and final portions of $l$ and the edges of $G$ as they appear around $x$, so we have at most 20 possibilities. (For the initial portion, we have 4 possibilities and for the final one 5.) For the sequence of edges, we have at most $|E(G)|^j$ possibilities. Taking into account that $0\le j\le 2k-2$ and $|E(G)|\le 6k-4$, the number of different signatures of the loops in $H\setminus L$ is smaller than $(6k)^{2k}$.
\smallskip

\begin{figure}[ht]
\begin{center}
\scalebox{0.5}{\includegraphics{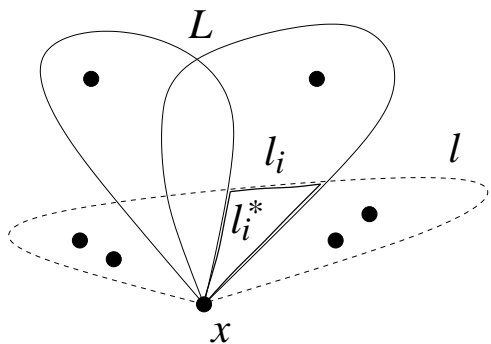}}
\caption{The definition of the $x$-loops $l^*_i$.}\label{leash}
\end{center}
\end{figure}

Next, we fix a signature and bound the number of elements in the subset $H^*\subseteq H\setminus L$ of all loops that have this signature. Any element $l\in H^*$ that has $j$ crossings with the edges of $G$, is divided into $j+1$ {\em curve-segments} (or, simply, {\em segments}) $l_0,l_1,\dots,l_j$. We extend each $l_i$ into an $x$-loop $l^*_i$. as follows. Let $l^*_i$ start with the standard path from $x$ to the initial point of $l_i$, followed by $l_i$, and then completed by the {\em reverse} of the standard path from $x$ to the final point $b$ of $l_i$. See Fig.~\ref{leash}. Note that the product $l_0^*l_1^*\ldots l_j^*$ is an $x$-loop homotopic to $l$. This means that for any two distinct (and, therefore, non-homotopic) loops $l,l'\in H^*$, there must be an index $0\le i\le j$ such that (the extension of) the $i$th segment of $l$ is not homotopic to (the extension of) the $i$th segment of $l'$.
\smallskip

We claim that for any fixed $i$, the number of distinct homotopy classes of the $i$th segments of the loops in $H^*$  is at most $f(n-1,k)$. If true, this would immediately imply that $|H^*|\le (f(n-1,k))^{j+1}$. Summing this bound over all signatures would eventually imply that
$$f(n,k)=|H|\le (6k)^{2k}(f(n-1,k))^{2k}.$$
\smallskip

It remains to prove the claim. We fix $i$ and a subset $H_0\subseteq H^*$ such that the $i$th segments of the loops in $H_0$ are pairwise non-homotopic. Let $F$ be the $L$-face (i.e., face of the drawing of the graph $G$) that contains the $i$th segment $l_i$ of a loop $l\in H^*$, and let $f$ denote the edge at which $l_i$ starts. Let us fix a point $x'\in F$ very close to $f$. (For $i=0$, the segment $l_i$ starts at $x$, between two edges of $G$, consecutive in the cyclic order. Then we pick $x'$ very close to $x$, between these two consecutive edges.) Assign to each $l\in H_0$ an $x'$-loop $l^*$ in $F\setminus T$, as follows. First, $l^*$ follows $f$ very closely till it reaches the $i$th curve-segment $l_i$ of $l$ close to its starting point. The second piece of $l^*$ follows $l_i$ almost to its endpoint, and then its third piece follows the boundary of $F$ very closely to get back to $p'$. If $l_i$ ends on the same edge $f$ of $G$ where it starts, the third piece of $l^*$ follows $f$ very closely. Otherwise, it follows the boundary of $F$ in a fixed cyclic direction.
\smallskip

If the pieces of $l^*$ that follow the boundary of $F$ run closer to it than the distance of any point in $T\cap F$ from the boundary of $F$, then the homotopy type of $l^*$ determines the homotopy type of $l_i^*$, and hence all the $|H_0|$ $x'$-loops will be pairwise non-homotopic. We can also choose these new loops in such a way  that every self-intersection of $l^*$ is also a self-intersection of $l_i$, and hence there are fewer than $k$ such self-intersections. In a similar manner, we can make sure that every intersection between two new loops is actually an intersection between the corresponding loops in $H_0$, and hence any two new loops intersect fewer than $k$ times. These new loops are pairwise non-homotopic in $S$. All of them lie in $F\setminus T\subset S$, therefore they are also non-homotopic there. Since $F$ is homeomorphic to the plane, $F\setminus T$ can be obtained from the plane by discarding $|F\cap T|$ points. We know that $|F\cap T|\le n-1$, because $F$ is an $L$-face and $L$ is balanced. This completes the proof of the claim and, hence, the lemma. \hfill  $\Box$
\medskip

Fix any $k\ge 1$. According to Corollary~\ref{bazis}, the upper bound in Theorem~\ref{diffloop} holds for $n=1$ and any $k\ge 1$. Let $n\ge 2$ and suppose that we have already verified the inequality $f(n-1,k)<2^{(2k)^{2(n-1)}}$. By Lemma~\ref{rekurzio}, we obtain
$$f(n,k)\le(6kf(n-1,k))^{2k}<(6k2^{(2k)^{2(n-1)}})^{2k}<2^{(2k)^{2n}},$$
completing the proof of Theorem~\ref{diffloop}.



\newpage

\section*{Appendix}

\subsection*{Proof of Lemma \ref{claim0}}

Assume for contradiction that there is a loose non-homotopic multigraph $H$, which is a counterexample to the lemma. We may assume the ${\rm cr}(H)$ is finite, as this can be achieved by infinitesimal perturbation. We choose $H$ to be a counterexample with {\em minimum crossing number}. If this number is zero (that is, the edges of $H$ have no self-intersections), we further minimize the number of {\em connected components }in $H$. Let $n$ stand for the number of vertices of $H$, the minimal counterexample, and let $m$ stand for the number of its edges.
\smallskip

Assume first that there is no self-intersecting edge in $H$, so we deal with a {\em planar} drawing.
In this case, we can also assume that $H$ is connected, otherwise two components could be
joined by an extra edge without creating a crossing.
(Note that this argument fails if we permit self-intersecting edges, as they may prevent
the addition of such an edge between two components without creating a crossing, see Fig. 1.)
Thus, the boundary of each face of $H$ can be visited by a single walk. These walks collectively cover every edge twice, so if each of them have at least three edges and the number of faces is $s$, then we have $3s\le2m$. Combining this inequality with Euler's formula $n+s=m+2$ gives $m\le3n-6$, which contradicts our assumption that $H$ was a counterexample. Therefore, $H$ must have a face bounded by a walk consisting of one or two edges. A boundary walk consisting of a single edge is a trivial loop, which is not permitted in a non-homotopic graph. A boundary walk of two edges is typically formed by two parallel edges that are homotopic, which is also disallowed in a non-homotopic graph. The only possibility is that the walk is back and forth along the same edge. In this case, we have $n=2$ and $m=1$, and $H$ is not a counterexample.
\smallskip

Therefore, our minimal counterexample $H$ must have at
least one self-inter\-secting edge $e$.
Find a minimal interval $\gamma$ of $e$ between two
occurrences of the same intersection point $p$.
This is a simple closed curve in the plane avoiding all vertices.
It partitions the sphere $S^2$ into two connected components.
We call them (arbitrarily) the \emph{left} and \emph{right sides}
of $\gamma$.
Obviously, $e$ is the only edge that may run between these sides.
Let $H_1$ and $H_2$ be the subgraphs of $H\setminus\{e\}$
induced by the vertices in the left and right sides of
$\gamma$, respectively. Both of them are loose topological multigraphs,
but they may contain homotopic edges.
By adding $p$ to both of them as an isolated vertex,
they become non-homotopic.
If an endpoint $u$ of $e$ lies in the left part,
then by adding to $H_1$ a non-self-intersecting edge
connecting $p$ and $u$ along $e$, we create no new
intersection and do not violate the non-homotopic condition either.
The resulting topological multigraph $H'_1$ is a loose
non-homotopic multigraph on the sphere with $n_1$ vertices and $m_1$ edges.
Analogously, we can construct the loose non-homotopic multigraph
$H'_2$ from $H_2$. Denote its number of vertices and edges by
$n_2$ and $m_2$, respectively. We have $n_1+n_2=n+2$
and $m_1+m_2\ge m$.
We eliminated a self-crossing (of $e$) and did not add any new crossings, so
the crossing numbers of both $H'_1$ and $H'_2$ are smaller than ${\rm cr}(H)$.
\smallskip

If $n_1, n_2>2$, then we have $m_1\le3n_1-6$ and $m_2\le3n_2-6$,
by the minimality of $H$. Summing up these inequalities, we get $m\le3n-6$, contradicting our assumption that $H$ was a counterexample.

If $n_1=1$ or $n_2=1$, all vertices of $H$ lie on the same side of
$\gamma$. In this case, by deleting $\gamma$ from $e$, the homotopy class of $e$ remains the same. Hence, the resulting topological multigraph is still a loose non-homotopic multigraph with $n$ vertices and $m$ edges, but its crossing number is smaller than that of $H$, contradicting the minimality of $H$.
\smallskip

Finally, consider the case $n_1=2$ or $n_2=2$.
By symmetry, we can assume that $n_1=2$, $n_2=n$, so we have a single vertex
$u$ of $H$ on the left side of $\gamma$ and $n-1$ vertices on the right side.
Note that no edge of $H\setminus\{e\}$ can lie in the left side. Indeed,
such an edge would be a trivial loop. If $e$ has at least one endpoint
in the right part, then we have $m_2=m$. This implies that $H'_2$ is another
counterexample to the lemma with fewer crossings,
contradicting the minimality of $H$.
\smallskip

Therefore, $e$ must be a loop at $u$.
The image of $e$ must separate a pair of vertices, $v, w\neq u$
of $H$ from each other,
as otherwise $e$ would be a trivial loop.
However, then we could draw another loop $e'$ along or very close to some parts of
$e$ with no self-intersection, so that it also separates $v$ and $w$.
Therefore, $e'$ is not trivial either.

Let $H'$ be the topological multigraph  obtained from $H$ by replacing $e$ by $e'$.
The loops $e$ and $e'$ are not necessarily homotopic,
but there is no other edge in $H'$ homotopic to $e'$, because there is no other loop at $u$.
Hence, $H'$ is a lose  non-homotopic
multigraph. This contradicts the minimality of $H$, because $H'$ has the same number of vertices and edges as $H$ does, but its crossing number is smaller.
This contradiction proves the lemma. \hfill  $\Box$

\subsection*{Proof  of Theorem~\ref{lowerbound2}}\label{main}


Let $x\in S^2$ and consider a family $L$ of $x$-loops in $S^2$ that start and end at $x$, but do not pass through $x$. With infinitesimal perturbations of the elements of $L$ and without creating any further intersection, one can attain that all intersections are simple: no point other than $x$ appears more than twice on the same loop or on different members of $L$. This will be assumed for all families of $x$-loops used in the rest of this section. A closed curve in $S^2$ is said to be an \emph{$L$-circle} if it is either a segment of a loop $l\in L$ between two appearances of a self-intersection point of $l$, or it consists of two segments of the same loop or two segments belonging to different loops in $L$, connecting the same pair of intersection points. If the two segments belong to the same loop, they are not allowed to overlap. We call a family of $L$-circles \emph{non-overlapping} if no two members of the family share a segment.

\begin{claim}\label{sokcircle}
Let $L$ be a family of $x$-loops consisting of a single loop with at least $k$ self-intersections or consisting of two loops intersecting each other at least $k$ times.

Then there is a non-overlapping family of $L$-circles, consisting of at least $k^{1/3}-1$ members.
\end{claim}

\noindent{\bf Proof.}
Suppose first that $L$ consists of two $x$-loops, $l_1$ and $l_2$. By the Erd\H os-Szekeres lemma~\cite{ErS35}, we can find $k'\ge\sqrt k$ intersection points $a_1, a_2, \ldots, a_{k'}$ that appear either in this order or in the reverse order on both $l_1$ and $l_2$. In this case, the segments of $l_1$ and $l_2$ between $a_i$ and $a_{i+1}$ form an $L$-circle, for each $1\le i<k'$. The family of these $L$-circles is non-overlapping, as claimed.
\smallskip

Alternatively, suppose that $L=\{l\}$ is a singleton family. The segment of $l$ between the two appearances of a self-intersection point $a$ is an $L$-circle. If there are at least $k^{1/3}$ among these $L$-circles that form a non-overlapping family, then we are done. If this is not the case, then we have a point $p\ne x$ on $l$ which is not an intersection point, but appears in at least $k^{2/3}$ of these single-part $L$-circles. That is, at least $k^{2/3}$ intersection points appear both on the initial segment $l_1$ of $l$ ending at $p$, and on the final segment $l_2$ of $l$, starting at $p$. We can then argue as we did in the case $|L|=2$. Using the Erd\H os-Szekeres lemma, we can find $k'\ge k^{1/3}$ intersection points $a_1,a_2,\dots, a_{k'}$ that appear either in this order or its reverse both on $l_1$ and $l_2$. For every $1\le i<k'$, the segments of $l_1$ and $l_2$ between $a_i$ and $a_{i+1}$ form an $L$-circle, and the family of these $L$-circles is non-overlapping. This finishes the proof. \hfill $\Box$

\begin{claim}\label{korbolmetszes}
Let $L$ be a family of $x$-loops not passing through a point $p\in S^2$, $p\ne x$. Let $L_1$ and $L_2$ be two disjoint sub-families of $L$ such that for $i\in\{1,2\}$ there is a non-overlapping family $C_i$ consisting of $k$ $L_i$-circles, each of which separates $x$ from $p$.

Then the total number of intersections between a loop in $L_1$ and a loop in $L_2$ is at least $k$.
\end{claim}

\noindent{\bf Proof.}
Let $F$ be the $(L_1\cup L_2)$-face which contains $p$. Let $q$ be an arbitrary non-intersection point on the boundary of $F$. Obviously, $q$ lies on an $x$-loop $l\in L_i$ with $i=1$ or $2$. Thus, $q$ belongs to the $L_{3-i}$-face containing $p$, and all $L_{3-i}$-circles in $C_{3-i}$ separate $x$ from $q$. The loop $l$ connects $x$ to $q$, so it must intersect all of these $L_{3-i}$-circles. As $C_{3-i}$ is a non-overlapping family, these intersections must be distinct intersection points of $l$ and some $x$-loop in $L_{3-i}$. This proves the claim. \hfill $\Box$
\medskip

Let us fix $n>1$ and a set $T$ of $n$ points on the $2$-sphere $S^2$. Let $S=S^2\setminus T$, and fix a point  $x\in S$.

\begin{lemma}\label{loopmetszes}
The minimal number $a_n(m)$ of crossings among $m$ pairwise non-homotopic $x$-loops in $S$ is super-quadratic in $m$. That is, for any fixed $n>1$, we have
 $$\lim_{m\rightarrow\infty}\frac{a_n(m)}{m^2}=\infty.$$
\end{lemma}

\noindent{\bf Proof.}
Let $L$ be a collection of $m$ non-homotopic $x$-loops in $S$ with the minimum overall number, $a_n(m)$, of crossings.

Choose the largest $k$ such that in any collection of at least $m/2$ non-homotopic $x$-loops in $S$,  there is one with at least $k$ self-crossings or two that cross each other at least $k$ times. By Theorem~\ref{diffloop}, $k$ goes to infinity as $m$ does.

We greedily divide the loops of $L$ into {\em blocks}, as follows. Each block is either a single loop crossing itself at least $k$ times, or a pair of loops crossing each other at least $k$ times. We do not use the same loop of $L$ twice.
Having formed at most $m/4$ blocks, we still have at least $m/2$ unused loops, and, by the definition of $k$, we can form yet another block. Therefore, the greedy procedure yields at least $m/4$ blocks.
\smallskip

By Claim~\ref{sokcircle}, for each block $B$, one can find a non-overlapping collection $C_B$ of at least $k^{1/3}-1$ $B$-circles. A $B$-circle $\gamma$ is called \emph{trivial} if it does not separate $x$ from any point of $T=S^2\setminus S$. The existence of a trivial $B$-circle would contradict the minimality of the total number of crossings in the collection $L$ of $x$-loops. Indeed, if a trivial $B$-circle consists of a single segment of an $x$-loop $l$, then deleting this segment does not affect the homotopy type of $l$, but decreases the number of crossings. If a trivial $B$-circle consists of two segments, then interchanging these segments does not affect the homotopy types of the corresponding loops. Now the number of crossings can be reduced by an infinitesimal perturbation of the original family $L$ or of the family obtained by this switch.
\smallskip

As no $B$-cycle in $C_B$ is trivial, for every block $B$, we can find a point $p_B\in T$ such that the number of $B$-cycles in $C_B$ which separate $p_B$ from $x$ is at least $k'=\frac{k^{1/3}-1}n$. There are only $n$ points in $T$, so there exists $p\in T$ such that $p_B=p$ for at least $m/(4n)$ blocks. By Claim~\ref{korbolmetszes}, any two distinct blocks $B$ and $B'$ for which $p_B=p_{B'}$, cross each other at least $k'$ times. This gives a total of at least $\binom{\lceil m/(4n)\rceil}2k'$ crossings. Since $n$ is fixed and $m$ tends to infinity, we know that $k$ and, therefore, $k'$ also tend to infinity. Thus, the number of crossings super-quadratic in $m$, as claimed. \hfill $\Box$
\medskip

\noindent{\bf Proof of Theorem~\ref{lowerbound2}.}
Let $G$ be a non-homotopic multigraph with $n$ vertices, $m$ edges, and with the smallest possible crossing number ${\rm cr}(n,m)$. As before, we can assume that there is no triple intersection among the edges, because we can get rid of these by infinitesimal perturbations. Obviously, we can find a set $E'$ of $m'\ge m/n^2$ parallel edges in $G$. We fix such a set $E'$, and in the rest of the proof we ignore all other edges of $G$. There are two cases.
\smallskip

{\em Case A:}  $E'$ consists of loops at a vertex $x$. In this case, choose a point $p$ very close to $x$ but not on any of the loops in $E'$. Let $S$ be the set obtained from the plane by deleting $p$ and all vertices of $G$ except $x$. The edges in $E'$ are pairwise non-homotopic $x$-loops in $S$. As $S$ can be obtained from the sphere $S^2$ by deleting $n+1$ points, these loops determine at least $a_{n+1}(m')$ intersections. According to Lemma~\ref{loopmetszes}, for a fixed $n$, this quantity is super-quadratic in $m'$ and, hence, also in $m$.
\smallskip

{\em Case B:} $E'$ consists of edges between two distinct vertices, $x$ and $y$. In this case, we pick two points, $p$ and $q$, very close to $x$ and $y$, respectively, which do not lie on any edge in $E'$. Now choose $S$ to be the set obtained from the plane by deleting all vertices of $G$ except $x$ and $y$, and also deleting $p$ and $q$. Any two edges of $E'$ form an $x$-loop. Moreover, for any $e_1,e_2,e_3\in E'$ with $e_2\ne e_3$, the $x$-loop formed by $e_1$ and $e_2$ is not homotopic in $S$ to the $x$-loop formed by $e_1$ and $e_3$.
\smallskip

We build a collection of pairwise non homotopic $x$-loops by pairing up edges of $E'$ in a greedy way, using every edge at most once. Suppose that the process stops with a collection $L$ of $m''$ $x$-loops. There are $m'-2m''$ unused edges left in $E'$. Fix any one of them, and combine it with each of the remaining ones to obtain $m'-2m''-1$ pairwise non-homotopic $x$-loops. Since we were unable to extend $L$ by another $x$-loop, each of these $x$-loops is homotopic to one of the $m''$ loops we have constructed so far. Therefore, we have $m''\ge m'-2m''-1$, and $m''\ge(m'-1)/3$.
\smallskip

All $x$-loops in $L$ pass through $y$. With an infinitesimal perturbation, one can get rid of this  multiple intersection without changing the homotopy classes of the $x$-loops or creating any additional intersection. Denote the resulting family of $x$-loops by $L'$. All loops in $L$ intersected at $y$. This may introduce up to $\binom{m''}2$ intersections between loops in $L'$ close to $y$. All other intersections among the members of $L'$ correspond to actual intersections between edges in $E'$.
\smallskip

Just like in Case A, $S$ can be obtained from the sphere by removing $n+1$ points. Hence, altogether there are at least $a_{n+1}(m'')$ intersections between the loops in $L'$, and the number of intersections between the edges of $E'$ is at least $a_{n+1}(m'')-\binom{m''}2$. In view of Lemma~\ref{loopmetszes}, this is super-quadratic in $m''$ and, hence, also in $m$. This completes the proof of Theorem~\ref{lowerbound2}. \hfill $\Box$

\end{document}